\documentclass[a4paper,12pt]{amsart}
\usepackage{bm}
\usepackage{ifthen}
\usepackage{graphicx}
\usepackage{mathrsfs}
\usepackage{color}
\nonstopmode \numberwithin{equation}{section}
\newtheorem{thm}{Theorem}[section]
\newtheorem{cor}[thm]{Corollary}
\newtheorem{lem}[thm]{Lemma}
\newtheorem{prop}[thm]{Proposition}

\newtheorem*{ThmA}{Theorem A}
\newtheorem*{ThmB}{Theorem B}

\theoremstyle{definition}

\newenvironment{rem}{%
\bigskip
\noindent \textsl{{\sl Remark. }}}{\bigskip}

\newenvironment{pf}[1][]{%
 \vskip 3mm
 \noindent
 \ifthenelse{\equal{#1}{}}%
  {{\slshape Proof. }}%
  {{\slshape #1.} }%
 }%
{\qed\bigskip}

\newcounter{alphabet}
\newcounter{tmp}

\newcommand{\A}{{\mathcal A}}

\newcommand{\C}{{\mathbb C}}

\newcommand{\D}{{\mathbb D}}

\newcommand{\R}{{\mathbb R}}

\newcommand{\es}{{\mathcal S}}

\newcommand{\U}{{\mathcal U}}

\newcommand{\tp}{{\null^\mathrm{T}}}

\newcommand{\aand}{{\quad\text{and}\quad}}

\newcounter{minutes}
\divide\time by 60
\newcounter{hours}
\multiply\time by 60 \addtocounter{minutes}{-\time}

\begin{document}
\bibliographystyle{amsplain}
\title[Sharp inequalities for logarithmic coefficients]{
Sharp inequalities for logarithmic coefficients and their applications
}


\author[S. Ponnusamy]{S. Ponnusamy}
\address{S. Ponnusamy, Department of Mathematics, Indian Institute of
Technology Madras, Chennai-600 036, India. }
\email{samy@iitm.ac.in}
\author[T. Sugawa]{Toshiyuki Sugawa}
\address{Graduate School of Information Sciences \\
Tohoku University \\
Aoba-ku, Sendai 980-8579, Japan}
\email{sugawa@math.is.tohoku.ac.jp}

\begin{abstract}
I.~M.~Milin proposed, in his 1971 paper, a system of inequalities for the logarithmic coefficients of normalized univalent functions
on the unit disk of the complex plane.
This is known as the Lebedev-Milin conjecture and
implies the Robertson conjecture which in turn implies the Bieberbach conjecture.
In 1984, Louis de Branges settled the long-standing Bieberbach conjecture by showing the Lebedev-Milin conjecture.
Recently, O.~Roth proved an interesting sharp inequality for the logarithmic coefficients based on the proof by de Branges.
In this paper, following Roth's ideas, we will show more general sharp inequalities
with convex sequences as weight functions and then establish several consequences of them.
We also consider the inequality with the help of de Branges system of linear
ODE for non-convex sequences where the proof is partly assisted by computer.
Also, we apply some of those inequalities to improve previously known results.
\end{abstract}
\thanks{
The present research was
supported by JSPS Grant-in-Aid for Scientific Research (B) 22340025 and JP17H02847.
The  work of the first author is supported by Mathematical Research Impact Centric
Support (MATRICS) of DST, India  (MTR/2017/000367).}
\maketitle

\section{Estimates of logarithmic coefficients}
Let $\A$ denote the set of normalized analytic functions on the open unit disk
$\D=\{z\in\C: |z|<1\}$ and $\es$ denote its subclass of univalent functions.
We define the logarithmic coefficients of $f$ by the formula
\begin{equation}\label{eq:log}
\log\frac{f(z)}{z}=2\sum_{n=1}^\infty \gamma_nz^n.
\end{equation}
Throughout the discussion, $\gamma _n:=\gamma _n(f)$ denote the logarithmic coefficients of a function $f\in \es$.
Louis de Branges \cite{dB85} solved the long-standing Bieberbach conjecture by
showing the Lebedev-Milin conjecture (see also \cite{FP85}):
For each $n\ge1,$
\begin{equation}\label{eq:LM}
\sum_{k=1}^n k(n-k+1)|\gamma_n|^2\le \sum_{k=1}^n \frac{n-k+1}k,
\end{equation}
where equality holds if and only if $f$ is the Koebe function $K(z)=z/(1-z)^2$
or its rotation $e^{-i\theta}K(e^{i\theta}z)=z/(1-e^{i\theta}z)^2$
for some $\theta\in\R.$
Note that for $f(z)=z/(1-e^{i\theta}z)^2$ we have
$\gamma_n=e^{in\theta}/n$ for $n=1,2,\dots.$

As an application of the de Branges theorem \eqref{eq:LM},
we will show a more general inequality.
As a preparation, we recall a notion of convexity for sequences.
A sequence of real numbers $p_n,~n=1,2,3,\dots,$ is called {\it convex}
if $p_n-2p_{n+1}+p_{n+2}\ge0$ for all $n\ge1.$
Note that $p_n=\varphi (n),~n=1,2,3,\dots,$ form a convex sequence if
$\varphi (x)$ is a convex function on $[1,+\infty)$ in the ordinary sense.
We can now state it as follows.

\begin{thm}\label{thm:simple}
Let $p_n,~n=1,2,3,\dots,$ be a convex sequence of non-negative numbers
with $p_1>0$ such that $\sum_{n=1}^\infty (p_n/n)<+\infty.$
For $f\in\es$ with expansion \eqref{eq:log}, the inequality
\begin{equation}\label{eq:main}
\sum_{n=1}^\infty np_n|\gamma_n|^2\le \sum_{n=1}^\infty\frac{p_n}{n}
\end{equation}
holds.
Moreover, the inequality is strict unless $f(z)$ has the form $z/(1-e^{i\theta}z)^2$
for some $\theta\in\R.$
\end{thm}

We remark that the theorem is not really new.
The same statement was already made by de Branges \cite{dB85} when the convex sequence
$p_n$ is eventually vanishing, i.e., $p_n=0$ for sufficiently large numbers $n.$
Zemyan in his 1993 paper \cite{Zem93} extended it to general convex sequences by approximating
them with eventually vanishing ones.
Therefore, he did not provide equality conditions.
For convenience of the reader, we give a direct proof of the theorem.

\begin{pf}[Proof of Theorem \ref{thm:simple}]
First note that $p_n/n\to0$ as $n\to\infty$ by the convergence assumption.
Put $q_n=p_n-p_{n+1}$ and
$$
\lambda_n=q_n-q_{n+1}=p_n-2p_{n+1}+p_{n+2}
$$
for $n=1,2,3,\dots.$
Then, by convexity, $\lambda_n\ge0$ and thus $q_n$ is a non-increasing sequence.
In particular, $q_n$ has a limit, say $q,$ as $n\to\infty.$
If $q\ne0,$ then $p_n$ is asymptotically equal to $nq,$ which
violates $p_n/n\to0.$
Hence, we conclude that $q=0.$
Since $q_n$ is non-increasing, we have
$p_n-p_{n+1}=q_n\ge0,$ which means $p_n$ is non-increasing.
In particular, $p_n$ has a limit, say $p,$ as $n\to\infty.$
Since $p_n\ge 0,$ we have $p\ge0.$
If $p>0,$ then $p_n\ge p,$ which implies $\sum p_n/n\ge
\sum p/n=+\infty,$ a contradiction.
Hence, the convergence assumption forces the sequence $p_n$ to converge to $0.$
Here we also note that, by the assumption $p_1>0,$
there is an $n\ge 1$ such that $\lambda_n>0.$

We now sum up the inequalities \eqref{eq:LM} with the weight $\lambda_n\ge0$
to obtain
\begin{equation}\label{eq:sum}
\sum_{n=1}^\infty \lambda_n\sum_{k=1}^n k(n-k+1)|\gamma_k|^2
\le \sum_{n=1}^\infty\lambda_n\sum_{k=1}^n \frac{n-k+1}k.
\end{equation}
Here, we note that
equality holds in \eqref{eq:sum} if and only if $f(z)=z/(1-e^{i\theta}z)^2,$
because equality must hold in \eqref{eq:LM} for at least one $n.$
The interchange of the order of summation gives us the inequality
$$
\sum_{k=1}^\infty k|\gamma_k|^2\sum_{n=k}^\infty \lambda_n(n-k+1)
\le\sum_{k=1}^\infty \frac1k\sum_{n=k}^\infty \lambda_n(n-k+1).
$$
If
\begin{equation}\label{eq:pk}
p_k=\sum_{n=k}^\infty\lambda_n(n-k+1),\quad k\ge1,
\end{equation}
then we would have the inequality \eqref{eq:main}.
We now show \eqref{eq:pk}.
Since $p_n\to0,$ we have
$$
p_k=\sum_{n=k}^\infty (p_n-p_{n+1})=\sum_{n=k}^\infty q_n.
$$
For convenience, for a fixed $k\ge1,$ we put $s_n=n-k$ for $n=k, k+1,\dots.$
Letting $N\ge k,$ we compute
\begin{align}\label{eq:1}
\notag
p_k-p_{N+1}&=\sum_{n=k}^N q_n=\sum_{n=k}^N(s_{n+1}-s_n)q_n \\
&=s_{N+1}q_N-s_kq_k+\sum_{n=k}^{N-1} s_{n+1}(q_n-q_{n+1}) \\
\notag
&=s_{N+1}q_N+\sum_{n=k}^{N-1} s_{n+1}\lambda_n.
\end{align}
Here, we used the fact that $s_k=0.$
In particular, we have
$$
s_{N+1}q_N+\sum_{n=k}^{N-1} s_{n+1}\lambda_n\le p_k.
$$
Since each term in the left-hand side is non-negative,
$$
\sum_{n=k}^\infty s_{n+1}\lambda_n
=\sup_{N>k}~\sum_{n=k}^{N-1} s_{n+1}\lambda_n
\le p_k <+\infty.
$$
Recalling $p_{N+1}\to0,$
we see by \eqref{eq:1} that $s_{N+1}q_N$ also has a limit,
say $b,$ as $N\to\infty.$
If $b\ne0,$ then $q_n$ is asymptotically $b/n$ and thus $p_n$ is
asymptotically $b\log n,$ which contradicts $p_n\to0.$
Thus we conclude that $b=0.$
Letting $N\to\infty$ in \eqref{eq:1}, we obtain the relation
$$
p_k=\sum_{n=k}^\infty s_{n+1}\lambda_n=\sum_{n=k}^\infty
(n-k+1)\lambda_n,
$$
and hence \eqref{eq:pk} is proved.
\end{pf}

In a recent paper by Roth \cite{Roth07}, he made the nice observation
that \eqref{eq:sum} could hold even if some of $\lambda_n$ are negative.
His idea is to show the inequality
\begin{equation}\label{eq:N}
\sum_{n=1}^N \lambda_n\sum_{k=1}^n k(n-k+1)|\gamma_k|^2
\le \sum_{n=1}^N\lambda_n\sum_{k=1}^n \frac{n-k+1}k
\end{equation}
for some $N\ge2$ by using the original idea of de Branges.
If $\lambda_n\ge0$ for $n>N,$ we obtain \eqref{eq:sum} by summing up
for $n>N$ with weight $\lambda_n.$
We will take a closer look at this case in the third section.


\section{Consequences of Theorem \ref{thm:simple}}

By various choices of positive convex sequences $p_n,$
we obtain many sharp inequalities on the logarithmic coefficients $\gamma_n$ of $f\in \es$.
The most fundamental one is perhaps $p_n=r^{2n}$ for a positive number $r.$
It is easy to check that this sequence satisfies the assumptions of Theorem \ref{thm:simple}
if and only if $r<1.$
Then we obtain the sharp inequality for the logarithmic area
$$
\sum_{n=1}^\infty n|\gamma_n|^2r^{2n}\le \sum_{n=1}^\infty\frac{r^{2n}}{n}
=\log\frac1{1-r^2},
$$
which is known as the Bazilevi\u c conjecture and proved by Milin and Grinshpan \cite{MG86}
(see also \cite{Milin91}).
The next fundamental example is $p_n=n^{-\alpha}$ for a constant $\alpha>0.$
Since $\varphi(x)=x^{-\alpha}$ is convex on $x>0,$ the sequence $p_n=\varphi(n)$
is convex. Therefore, as a corollary of Theorem \ref{thm:simple}, we obtain the inequality
$$
\sum_{n=1}^\infty n^{1-\alpha}|\gamma_n|^2\le
\sum_{n=1}^\infty\frac{1}{n^{1+\alpha}}=\zeta(\alpha+1),
$$
where $\zeta(x)$ denotes the Riemann zeta function.
Equality holds if and only if $f$ is a rotation of the Koebe function $z/(1-z)^2.$
This inequality was proved by Zemyan \cite[Theorem 3 (b)]{Zem93}.
Letting $\alpha=1$ in particular, we obtain the Duren-Leung inequality \cite{DL79}
\begin{equation}\label{eq:DL}
\sum_{n=1}^\infty |\gamma_n|^2\le \frac{\pi^2}6.
\end{equation}
It is worth recalling that this inequality was proved even before
de Branges' proof of the Lebedev-Milin conjecture.

We summarize other choices in the following lemma.

\begin{lem}\label{lem1}
For each choice of the following, the sequence $p_n~(n=1,2,3,\dots)$ is positive and convex.
\begin{enumerate}
\item[(1)] $p_n=\dfrac1{n+\alpha}$ and $\alpha>-1,$
\vspace{1mm}
\item[(2)] $p_n=\dfrac{n}{n^2+a n+b}$ for $a, b\in\R$
with $a+b+1>0,~a+3\ge0$ and $ (6+a)b\le 6.$
\vspace{1mm}
\item[(3)] $p_n=\dfrac{n}{(n+\alpha)(n+\beta)}$ for
$\alpha>-1,~\beta>-1$ with $(\alpha+\beta+6)\alpha\beta\le 6$ and $\alpha\beta\le 6.$
\vspace{1mm}
\item[(4)] $p_n=\dfrac{1}{n^2+a n+b}$ for $a, b\in\R$ with
$a+b+1>0,~ a+2\ge 0$ and $b\le a^2+6a+11.$
\vspace{1mm}
\item[(5)] $p_n=\dfrac{1}{(n+\alpha)(n+\beta)}$ for $\alpha>-1,~\beta>-1.$
\vspace{1mm}
\item[(6)] $p_n=\dfrac{n^2}{(n+\alpha)^2(n+\beta)}$ for $\alpha>-1,~\beta>-1$
with $|\alpha|(1+3|\beta|+\beta^2)\le 1/2.$
\vspace{1mm}
\item[(7)] $p_n=(n+\alpha)r^n$ for $\alpha>-1$ and $r\in(0,1)$ with
$2\le(\alpha+1)\log(1/r).$
\vspace{1mm}
\end{enumerate}
\end{lem}
\begin{pf}
We will take the following strategy to show the assertion.
First we choose a smooth function $\varphi$ so that $p_n=\varphi(n).$
If we confirm that $\varphi(x)$ is convex on $N\le x<+\infty$ for an integer $N\ge1,$
then it is enough
to check the condition $\lambda_n=p_n-2p_{n+1}+p_{n+2}\ge0$
for $n=1,2,\dots, N-1.$

(1) Since $\varphi(x)=1/(x+\alpha)$ is convex on $1\le x$ for $\alpha>-1,$
the assertion follows.

(2) First note that $p_n>0$ for $n\ge1$ by the first two conditions on parameters.
Indeed, $n^2+an+b=(n-1)(n+1+a)+1+a+b\ge 1+a+b>0$ for $n\ge 1.$
As a necessary condition, we have
$$
\lambda_1=\frac{2(6-6b-ab)}{(1+a+b)(4+2a+b)(9+3a+b)}\ge0,
$$
which is certainly implied by the assumption.
Let $\varphi(x)=x/(x^2+a x+b)$ and compute
$$
\varphi''(x)=\frac{2(x^3-3b x-ab)}{(x^2+ax+b)^3}.
$$
We note here that $b\le 6/(6+a)\le 2$ by the assumptions $a+3\ge0$
and $(6+a)b\le 6.$
Since
$$
x^3-3bx-ab=x^3-3bx+6b-6+(6-6b-ab)
\ge x^3-3bx+6b-6=:h(x),
$$
it is enough to show that $h(x)\ge 0$ for $x\ge 2.$
Since $h'(x)=3(x^2-b),$ the function $h(x)$ is increasing in $2\le x<+\infty$
and thus $h(x)\ge h(2)=2>0$ as required.

(3) We apply the previous case for $a=\alpha+\beta$ and $b=\alpha\beta$
to get the assertion.

(4) As in the case (2), we see that $p_n>0$ by the first two conditions on $a,b.$
Also, the inequality
$$
\lambda_1=\frac{2(a^2+6a+11-b)}{(1+a+b)(4+2a+b)(9+3a+b)}\ge0
$$
holds by assumption.
Let $\varphi(x)=1/(x^2+a x+b)$ and compute
$$
\varphi''(x)=\frac{2(3x^2+3ax+a^2-b)}{x^2+ax+b)^3}.
$$
Since $h(x)=3x^2+3ax+a^2-b$ is increasing in $x\ge -a/2~(\le 1),$ we obtain
$h(x)\ge h(2)=a^2+6a+12-b\ge 1>0$ for $x\ge 2.$
Thus we conclude that $\varphi(x)$ is convex on $2\le x<+\infty.$

(5) Just apply (4) with $a=\alpha+\beta$ and $b=\alpha\beta.$

(6) Let $\varphi(x)=x^2/[(x+\alpha)^2(x+\beta)].$
Then
$$
\varphi'(x)=\frac{2(x^4-2\alpha x^3-6\alpha\beta x^2-2\alpha\beta^2
+\alpha^2\beta^2}{(x+\alpha)^4(x+\beta)^3}.
$$
For $x\ge 1,$ we have
$$
x^4-2\alpha x^3-6\alpha\beta x^2-2\alpha\beta^2+\alpha^2\beta^2
\ge x^4(1-2|\alpha|-6|\alpha\beta|-2|\alpha\beta^2|)\ge0,
$$
which implies that $\varphi(x)$ is convex on $1\le x<+\infty.$

(7) It is enough to observe the formula $\varphi''(x)=\{2+(\alpha+x)\}r^x\log r$
for $\varphi(x)=(x+\alpha)r^x.$
\end{pf}

\begin{cor}\label{cor1}
For the logarithmic coefficients $\gamma_n$ of $f\in\es,$
the following inequalities hold. Each of them is strict unless
$f$ is not a rotation of the Koebe function $z/(1-z)^2.$
\begin{enumerate}
\item $\displaystyle \sum_{n=1}^\infty \frac{n}{n+\alpha }|\gamma_n|^2\le
\sum_{n=1}^\infty\frac{1}{n(n+\alpha)}=:A_\alpha$
for $\alpha>-1.$
When $\alpha\ne0,$ we have the expressions
$$
A_\alpha
=\frac{1}{\alpha}\int_0^1\frac{1-t^{\alpha}}{1-t}\, dt
=\frac{\psi(\alpha+1)-\psi(1)}{\alpha}.
$$
Here and in the sequel $\psi(x)=\Gamma'(x)/\Gamma(x)$ denotes the Digamma function.
In particular, letting $\alpha =1,2,3,1/2, -1/2,$ the following sharp inequalities are
deduced:
\begin{enumerate}
\renewcommand{\labelenumii}{[{\alph{enumii}}]}
\item$\displaystyle
~\sum_{n=1}^\infty \frac{n}{n+1}|\gamma_n|^2\le 1,$
\item$\displaystyle
~\sum_{n=1}^\infty \frac{n}{n+2}|\gamma_n|^2\le \frac{3}{4},$
\item$\displaystyle
~\sum_{n=1}^\infty \frac{n}{n+3}|\gamma_n|^2\le \frac{11}{18},$
\item$\displaystyle
~\sum_{n=1}^\infty \frac{n}{2n+1}|\gamma_n|^2\le 2(1-\log 2),$
\item$\displaystyle
~\sum_{n=1}^\infty \frac{n}{2n-1}|\gamma_n|^2\le 2\log 2.$
\end{enumerate}

\item $\displaystyle
\sum_{n=1}^\infty \frac{n^2}{n^2+\alpha ^2}|\gamma_n|^2\le
B_\alpha:=\frac{\pi\alpha \coth\pi\alpha -1}{2\alpha^2}$ for $\alpha \in (0,1].$

\item $\displaystyle
\sum_{n=1}^\infty \frac{n^2}{(n+\alpha)(n+\beta)} |\gamma_n|^2
\le C_{\alpha,\beta}$ for $\alpha,\beta\in(-1,+\infty)$
with $(\alpha+\beta+6)\alpha\beta\le 6$ and $\alpha\beta\le 6.$
Here,
$$
C_{\alpha,\beta} = \begin{cases}
\displaystyle \frac{1}{\beta-\alpha}
\int_0^1\frac{t^{\alpha}-t^{\beta}}{1-t}\,dt
=\frac{\psi(1+\beta)-\psi(1+\alpha)}{\beta-\alpha}
&\quad \text{if}~ \alpha\neq \beta
\vspace{2mm}
\\
\displaystyle \sum_{n=1}^\infty\frac1{(n+\alpha)^2}
=\int_0^1\frac{t^{\alpha}\log (1/t)}{1-t}\,dt=\psi'(1+\alpha),
&\quad \text{if}~ \alpha=\beta.
\end{cases}
$$
In particular,
\begin{enumerate}
\renewcommand{\labelenumii}{[{\alph{enumii}}]}
\item$\displaystyle
~\sum_{n=1}^\infty \frac{n^2}{4n^2-1}|\gamma_n|^2\le \frac{1}{2},
$
\item$\displaystyle
~\sum_{n=1}^\infty \frac{n^2}{(n+1)(2n+1)}|\gamma_n|^2\le 2\log2-1.
$
\end{enumerate}

\item $\displaystyle
\sum_{n=1}^\infty \frac{n}{(n+\alpha)(n+\beta)} |\gamma_n|^2\le D_{\alpha,\beta}$
for $\alpha,\beta\in(-1,+\infty),$ where
$$
D_{\alpha,\beta}=
\int_0^1\frac{\beta(1-t^{\alpha}) -\alpha(1-t^{\beta})}%
{\alpha \beta(\beta-\alpha)(1-t)}\,dt
=-\frac1{\beta-\alpha}\left(\frac{\psi(1+\beta)}\beta-\frac{\psi(1+\alpha)}\alpha\right)
-\frac{\psi(1)}{\alpha\beta}
$$
for nonzero $\alpha,\beta$ with $\alpha\ne\beta,$
$$
D_{\alpha,0}=D_{0,\alpha}
=\frac{\zeta(2)-A_\alpha}{\alpha}
=\frac{\pi^2/6-A_\alpha}{\alpha}
$$
and
$$
D_{\alpha,\alpha}
= \int_0^1\frac{1-t^{\alpha}-\alpha t^{\alpha}\log (1/t)}{\alpha^2(1-t)}\,dt
=\frac{\psi(1+\alpha)-\psi(1)}{\alpha^2}-\frac{\psi'(1+\alpha)}{\alpha}
$$
for nonzero $\alpha,$ and
$$
D_{0,0}=\zeta(3).
$$
In particular,
\begin{enumerate}
\renewcommand{\labelenumii}{[{\alph{enumii}}]}
\item$\displaystyle
~\sum_{n=1}^\infty \frac{n}{(n+1)(n+2)}|\gamma_n|^2\le \frac14.$
\item$\displaystyle
~\sum_{n=1}^\infty \frac{1}{n+1}|\gamma_n|^2\le \frac{\pi^2}{6}-1,$
\item$\displaystyle
~\sum_{n=1}^\infty \frac{n}{(n+1)^2}|\gamma_n|^2\le 2-\frac{\pi^2}{6}.$
\end{enumerate}
\item $\displaystyle
\sum_{n=1}^\infty \frac{n^3}{(n+\alpha)^2(n+\beta)} |\gamma_n|^2
\le E_{\alpha,\beta}$ for $\alpha,\beta\ne0$ with $|\alpha|(1+3|\beta|+\beta^2)\le 1/2,$
where
$$
E_{\alpha,\beta}=\begin{cases}\displaystyle
\frac{\beta C_{\alpha,\beta}-\alpha C_{\alpha,\alpha}}{\beta-\alpha}
&\quad\text{if}~\alpha\ne\beta,
\vspace{2mm}
\\ \displaystyle
\psi'(1+\alpha)+\frac\alpha 2\psi''(1+\alpha) &\quad\text{if}~ \alpha=\beta.
\end{cases}
$$
\item$\displaystyle
\sum_{n=1}^\infty n(n+\alpha)|\gamma_n|^2r^{2n}
\le \frac{r^2}{1-r^2}+\alpha\log\frac1{1-r^2}$
for $\alpha>-1$ and $0<r<1$ with $1\le(\alpha+1)\log(1/r).$
\end{enumerate}
\end{cor}

\begin{pf} Basically, all the inequalities follow from Theorem \ref{thm:simple}
and Lemma \ref{lem1}.
The remaining task is only to compute the sum $\sum_{n=1}^\infty p_n/n.$

(1) By the formula
$$
A_\alpha=\sum_{n=1}^\infty\frac1{n(n+\alpha)}
=\frac1\alpha\sum_{n=1}^\infty \left(\frac1n-\frac1{n+\alpha}\right)
=\frac1\alpha\sum_{n=1}^\infty \int_0^1 \big(t^{n-1}-t^{n+\alpha-1}\big)dt
$$
we easily obtain the first expression.
The second expression can be obtained by the well-known formula (see \cite[6.3.16]{AS:hand})
$$
\psi(1+x)
=-\gamma+\sum_{n=1}^\infty \frac{x}{n(n+x)}
=-\gamma+\sum_{n=1}^\infty \left(\frac1n-\frac1{n+x}\right)
\quad(x\ne -1,-2,-3,\cdots),
$$
where $\gamma$ is Euler's constant.
The following formulae are convenient in practical computations:
$$
\psi(1+x)=\psi(x)+\frac1x
\aand
\psi(1)=-\gamma.
$$

(2) We need to show the identity
$$
\sum_{n=1}^\infty \frac{1}{n^2+\alpha ^2}
=  \frac{\pi\alpha \coth\pi\alpha -1}{2\alpha^2}.
$$
This can be deduced by subsituting $z=i\alpha$ into the well-known formula
(see \cite[p.~189]{Ahlfors:ca})
$$
\pi\cot \pi z=\frac1z+\sum_{n=1}^\infty\frac{2z}{z^2-n^2}.
$$

(3) The required formula
$$
C_{\alpha,\beta}=\sum_{n=1}^\infty\frac{1}{(n+\alpha)(n+\beta)}
$$
can be shown in the same way as in (2).
The particular cases follow from the computations
$C_{1/2,-1/2}=2$ and $C_{1/2,1}=2(2\log 2-1).$

(4) We need to check the formula
$$
D_{\alpha,\beta}=\sum_{n=1}^\infty \frac{1}{n(n+\alpha)(n+\beta)}.
$$
For the generic case $\alpha\neq\beta$, we may write the right-hand side  in the form
$$
\frac{1}{\beta-\alpha}\left (\sum_{n=1}^\infty \frac{1}{n(n+\alpha)}
 -\sum_{n=1}^\infty \frac{1}{n(n+\beta)} \right )
=-\left(\frac{A_\beta-A_\alpha}{\beta-\alpha}\right)
$$
and the assertion follows immediately from Case (2).
The rest of the assertions follows easily from a standard limiting process.

(5) We have only to use the expression
$$
\frac{n}{(n+\alpha)^2(n+\beta)}=
\frac{\beta}{\beta-\alpha}\frac1{(n+\alpha)(n+\beta)}
-\frac{\alpha}{\beta-\alpha}\frac1{(n+\alpha)^2}
$$
for $\alpha\ne\beta.$
The case when $\alpha=\beta$ follows from a suitable limiting process.

(6) Apply Lemma \ref{lem1} (7) with $r^2$ instead of $r.$ It is easy to check the formula
$$
\sum_{n=1}^\infty \frac{n+\alpha}{n}r^{2n}=\frac1{1-r^2}+\alpha\log\frac1{1-r^2}.
$$
\end{pf}

It is noteworthy that the above formulae of various series in the proof of the corollary
are valid in general regardless of the parameter conditions.

We remark that
$$
A_0=\int_0^1\frac{\log(1/t)}{1-t}dt=\psi'(1)=\sum_{n=1}^\infty\frac1{n^2}
=\zeta(2)=\frac{\pi^2}{6}.
$$
Therefore, we have the Duren-Leung inequality \eqref{eq:DL} as the
limiting case as $\alpha\to0$ in (2).
Also, we should confess that an application of Lemma \ref{lem1} (4)
could not be included in the corollary due to difficulty of evaluation of infinite
series of the form
$$
\sum_{n=1}^\infty\frac1{n(n^2+an+b)}
$$
when $a^2-4b<0.$
We add a couple of further consequences of Theorem \ref{thm:simple}.

\begin{cor}
\begin{enumerate}
\item$\displaystyle
\sum_{n=1}^\infty \frac{n^2}{(n+1)^3}|\gamma_n|^2\le \zeta(3)-1,$
\item$\displaystyle
\sum_{n=1}^\infty \frac{n}{(n+1)^3}|\gamma_n|^2\le \frac{1}{6}\left(18-\pi^2-6\zeta(3)\right).$
\end{enumerate}
\end{cor}

\begin{pf}
(1) follows from the fact that $\varphi(x)=x/(x+1)^3$ is convex on $1\le x<+\infty.$
(2) follows also from the convexity of $\varphi(x)=1/(x+1)^3$ and the computation
$$
\sum_{n=1}^\infty \frac1{n(n+1)^3}
=\sum_{n=1}^\infty\left(\frac{1}{n(n+1)^2}-\frac1{(n+1)^3}\right)
=C_{1,1}-(\zeta(3)-1).
$$
\end{pf}

\section{Computer-assisted proof of the inequality for non-convex sequences}

In the first section, we presented an inequality of the logarithmic
coefficients $\gamma_n$ for a convex sequence $p_n.$
The inequality may hold even if $p_n$ is not convex; namely,
some of $\lambda_n=p_n-2p_{n+1}+p_{n+2}$ are negative.
We review the idea due to Roth \cite{Roth07} and then reformulate it
in a convenient form so that one can check the conditions by using computers.

We recall the proof of the Lebedev-Milin conjecture \eqref{eq:LM}
by following FitzGerald and Pommerenke \cite{FP85}.
Fix $n\ge1.$
The key idea is to consider the de Branges system of linear ODE:
$$
\tau_{n,k}(t)-\tau_{n,k+1}(t)
=-\frac{\tau_{n,k}'(t)}k-\frac{\tau_{n,k+1}'(t)}{k+1},
\quad \tau_{n,k}(0)=n-k+1
$$
for $k=1,2,\dots, n,$ where we put $\tau_{n,n+1}(t)\equiv 0.$
With the aid of L\"owner chains, we can see that \eqref{eq:LM} follows
from the inequalities $\tau_{n,k}'(t)<0,~ t>0, k=1,2,\dots, n.$
See \cite{FP85} for details.
It is known that $\tau_{n,k}'(t)$ can be expressed in terms of
Jacobi polynomials (see \cite[(2.3)]{FP85}):
\begin{equation}\label{eq:J}
\tau_{n,k}'(t)=-ke^{-kt}\sum_{j=0}^{n-k}P_j^{(2k,0)}(1-2e^{-t}).
\end{equation}
Here, Jacobi polynomials are defined, for instance, by Rodrigues' formula
$$
P_j^{(\alpha,\beta)}(x)=\frac{(-1)^n}{2^nn!}(1-x)^{-\alpha}(1+x)^{-\beta}
\frac{d}{dx}\left[(1-x)^\alpha(1+x)^\beta(1-x^2)^n\right].
$$
The Askey-Gasper inequality was a key step to confirm
$\tau_{n,k}'(t)<0.$

Roth \cite{Roth07} observed that the same idea works for the inequality
\eqref{eq:N}.
Namely, consider the solution to the initial value problem
\begin{equation}\label{eq:de}
\tau_{k}(t)-\tau_{k+1}(t)
=-\frac{\tau_{k}'(t)}k-\frac{\tau_{k+1}'(t)}{k+1},
\quad \tau_{k}(0)=\mu_k
\end{equation}
for $k=1,2,\dots, N,$
where $\tau_{N+1}=0$ and
$$
\mu_k=\sum_{j=1}^{N-k+1}j\lambda_{j+k-1}=\sum_{n=k}^N\lambda_n(n-k+1).
$$
If the condition
\begin{equation}\label{eq:ini}
\tau_k'(t)<0~ \text{for}~ t>0,~ k=1,2,\dots, N,
\end{equation}
holds, then \eqref{eq:N} can be
deduced in the same way as in \cite{FP85} (see \cite{Roth07} for details).
When $p_n=n/(n+1)^2$ and $\lambda_n=p_n-2p_{n+1}+p_{n+2},$
by solving the differential equations, Roth \cite{Roth07} showed that the
condition \eqref{eq:ini} holds for $N=5.$

We take now a slightly different approach below.
In view of the form of \eqref{eq:N},  we see that $\tau_k$
can be described in terms of the original $\tau_{n,k}$'s.
Indeed, we have
\begin{equation}\label{eq:rep}
\tau_k=\sum_{n=k}^N \lambda_n\tau_{n,k}.
\end{equation}
Therefore, by \eqref{eq:J},
$\tau_k'$ can be expressed in terms of Jacobi polynomials:
$$
\tau_k'(t)
=- ke^{-kt}\sum_{n=k}^N \lambda_n
\sum_{j=0}^{n-k}P_j^{(2k,0)}(1-2e^{-t})
=- ke^{-kt}\sum_{j=0}^{N-k}\nu_{k,j} P_j^{(2k,0)}(1-2e^{-t}),
$$
where
\begin{equation*}
\nu_{k,j}=\sum_{n=j+k}^{N}\lambda_n=\lambda_{j+k}+\lambda_{j+k+1}+\dots+
\lambda_N =q_{j+k}-q_{N+1}.
\end{equation*}
We can now summarize these observations as the following theorem.

\begin{thm}\label{thm:N}
Let $p_n,~n=1,2,3,\dots,$ be a sequence of non-negative numbers and
set $q_n=p_n-p_{n+1}$ and $\lambda_n=q_n-q_{n+1}.$
Suppose that there exists a number $N\ge1$ satisfying the following three  conditions:
\begin{enumerate}
\item[(0)] $p_{N+1}>0,$
\item[(i)] $\lambda_n\ge0$ for $n>N,$
\item[(ii)] \label{item:ii}
$\displaystyle Q_k(x)=\sum_{j=0}^{N-k}\nu_{j+k} P_j^{(2k,0)}(x)>0$
for $-1<x<1$ and $k=1,2,\dots, N,$
where $\nu_m=q_m-q_{N+1}.$
\end{enumerate}
Then the inequality
$$ \sum_{n=1}^\infty np_n|\gamma_n|^2\le \sum_{n=1}^\infty\frac{p_n}{n}
$$
holds. Here, equality holds precisely when $f$ is a rotation of the Koebe function $z/(1-z)^2.$
\end{thm}

As an example, let us look at the case of Roth \cite{Roth07}.
Let $p_n=n/(n+1)^2.$
Then $\lambda_1=-1/144<0$ but $\lambda_n>0$ for $n>1.$
Take $N=5$ and compute $Q_k,~k=1,2,3,4,5,$ as follows:
\begin{align*}
Q_1(x)&=\dfrac{153191+313428x+443802x^2+517076x^3+249375x^4}{5644800} \\
Q_2(x)&=\dfrac{38929+77359x+82447x^2+35625x^3}{705600} \\
Q_3(x)&=\dfrac{139643+218986x+106875x^2}{2822400} \\
Q_4(x)&=\dfrac{15171+11875x}{705600} \\
Q_5(x)&=\dfrac{95}{28224}.
\end{align*}
By numerical computations, we can check that $Q_k(x)$ has no roots
on the interval $(-1,1).$
Hence, we verified the Roth inequality \cite{Roth07}
\begin{equation}\label{eq:Roth}
\sum_{n=1}^\infty \frac{n^2}{(n+1)^2}|\gamma_n|^2\le
\frac{\pi^2}{6}-1.
\end{equation}

It is not necessarily easy to check condition (ii) in the theorem.
Indeed, we have no general idea about how large $N$ should be chosen.
Therefore, the following necessary condition is useful in practical tests.

\begin{prop}\label{prop:nec}
Under the hypothesis of Theorem \ref{thm:N},
a necessary condition for \eqref{item:ii} is
$$v_k=v_{k,N}=\sum_{j=0}^{[(N-k)/2]}\lambda_{k+2j}
=\lambda_k+\lambda_{k+2}+\lambda_{k+4}+\dots+\lambda_{N'}\ge0,
\quad k=1,2,\dots, N,
$$
where $N'=N$ if $N-k$ is even and $N'=N-1$ if $N-k$ is odd.
\end{prop}
\begin{pf}
For \eqref{item:ii}, the condition $\tau_k'(0)\le0$ is necessary.
It is noted in \cite[p.~685]{FP85} that $\tau_{n,k}'(0)=-k$
if $n-k$ is even and $\tau_{n,k}'(0)=0$ if $n-k$ is odd.
By \eqref{eq:rep}, we have
$$
\tau_k'(0)=\sum_{n=k}^N\lambda_n\tau_{n,k}'(0)
=-k(\lambda_k+\lambda_{k+2}+\cdots+\lambda_{N'})=-kv_k.
$$
Thus we have the condition $v_k\ge0.$
\end{pf}

For instance,
\begin{align*}
v_{1,1}&=\lambda_1, \\
v_{1,2}&=\lambda_1, ~ v_{2,2}=\lambda_2, \\
v_{1,3}&=\lambda_1+\lambda_3, ~ v_{2,3}=\lambda_2, ~ v_{3,3}=\lambda_3, \\
v_{1,4}&=\lambda_1+\lambda_3, ~ v_{2,4}=\lambda_2+\lambda_4, ~ v_{3,4}=\lambda_3,~ v_{4,4}=\lambda_4.
\end{align*}

In particular, we observe that the choice $N\le 2$ does not work for Theorem \ref{thm:N}
when $\lambda_1<0.$

\begin{rem}
Unfortunately, the condition in the above proposition is not necessarily
sufficient.
Let $x_k(t)=\tau_k(t)/k$ for $k=1,2,\dots, N.$
Then the system of ODE \eqref{eq:de} turns to
$$
x_k'=-kx_k+2(k+1)x_{k+1}-2(k+2)x_{k+2}+\dots+(-1)^{N-k+1}2Nx_{N},\quad
x_k(0)=\mu_k/k
$$
for $k=1,2,\dots, N.$
Introducing the column vector $\bm{x}=(x_1,\dots,x_N)^{\mathrm T},$
the system can be expressed by $\bm{x}'=A_N\bm{x}$ for the $N\times N$
matrix $A_N$ corresponding to the above equation.
For example,
$$
A_2=\begin{pmatrix} -1 & 4 \\ 0 & -2 \end{pmatrix},\quad
A_3=\begin{pmatrix} -1 & 4 & -6 \\ 0 & -2 & 6 \\ 0 & 0 & -3\end{pmatrix}.
$$
Letting $\bm{x}_0$ be the initial vector at $t=0,$
the solution can be given by $\bm{x}=e^{tA_N}\bm{x}_0$
and thus $\bm{x}'=e^{tA_N}A_N\bm{x}_0.$
In our case, $A_N\bm{x}_0=-\tp(v_1,\dots, v_N),$ where $v_k$ are as in
Proposition \ref{prop:nec}.
Simple computations give us
$$
e^{tA_2}
=\begin{pmatrix} e^{-t} & 4e^{-t}(1-e^{-t}) \\ 0 & e^{-2t} \end{pmatrix},\quad
e^{tA_3}
=\begin{pmatrix} e^{-t} & 4e^{-t}(1-e^{-t}) & 3e^{-t}(1-e^{-t})(3-5e^{-t}) \\%
 0 & e^{-2t} & 6e^{-2t}(1-e^{-t}) \\ 0 & 0 & e^{-3t}\end{pmatrix}.
$$
We observe that the entry $3e^{-t}(1-e^{-t})(3-5e^{-t})$ of $e^{tA_3}$
takes negative values when $t>0$ is small enough.
If $\bm{v}=(v_1,v_2,v_3)^{\mathrm{T}}$ is very close to $(0,0,1)^{\mathrm{T}},$
then the first entry of $\bm{x}'(t)=e^{tA_3}A_3\bm{v}$ will take
negative values even if $v_k>0$ are satisfied.
\end{rem}

As an example, we consider the sequence
$$
p_n=\frac{n}{n^2+\alpha^2}\quad (n=1,2,3,\dots)
$$
for $\alpha>0,$ which appears in Corollary \ref{cor1} (2).
Put $q_n=p_n-p_{n+1}$ and $\lambda_n=q_n-q_{n+1}$ as before.
By Lemma \ref{lem1} and its proof, we see that the sequence $p_n$ is
convex if and only if $\alpha\le 1.$
It might be an interesting problem to find the largest value $\alpha$ so that
the inequality
\begin{equation}\label{eq:alpha}
\sum_{n=1}^\infty \frac{n^2}{n^2+\alpha ^2}|\gamma_n|^2
\le
\sum_{n=1}^\infty\frac{1}{n^2+\alpha^2}
=\frac{\pi\alpha \coth\pi\alpha -1}{2\alpha^2}
\end{equation}
holds for the logarithmic coefficients $\gamma_n$ of every function $f\in\es.$
For simplicity, put $b=\alpha^2>1.$
Then $\lambda_1<0$ but $\lambda_2\ge 0$ if $b\le 8/3.$
As we saw, we should choose $N\ge 3.$
When $N=3,$ we compute
$$
v_{1,3}=\lambda_1+\lambda_3
=\frac{12(440-317b-40b^2-3b^3)}{(1+b)(2^2+b)(3^2+b)(4^2+b)(5^2+b)}.
$$
Therefore, $b\le b_0$ is necessary and sufficient for $v_{1,3}\ge0,$
where $b_0=1.1925184\cdots$ is the unique real solution to the equation
$440-317b-40b^2-3b^3=0.$
In this case, $v_{2,3}=\lambda_2>0,~ v_{3,3}=\lambda_3>0.$
A numerical computation tells us that the polynomial $Q_1(x)$ in Theorem \ref{thm:N}
with $N=3$ and $b=b_0$ assumes a negative value on $-1<x<1,$ see Figure 1.
Therefore, the condition in Proposition \ref{prop:nec} is, indeed, not sufficient
for condition (ii) to hold in the theorem.
On the other hand, numerical experiments suggest that
$Q_1(x)>0$ on $-1<x<1$ for $b\le 1.19245.$
Other conditions $Q_2(x)>0$ and  $Q_3(x)>0$ can be checked more easily.
Thus, in this case, the inequality \eqref{eq:alpha} holds for $b\le 1.19245.$

\begin{figure}
\begin{center}
\includegraphics[width=100mm]{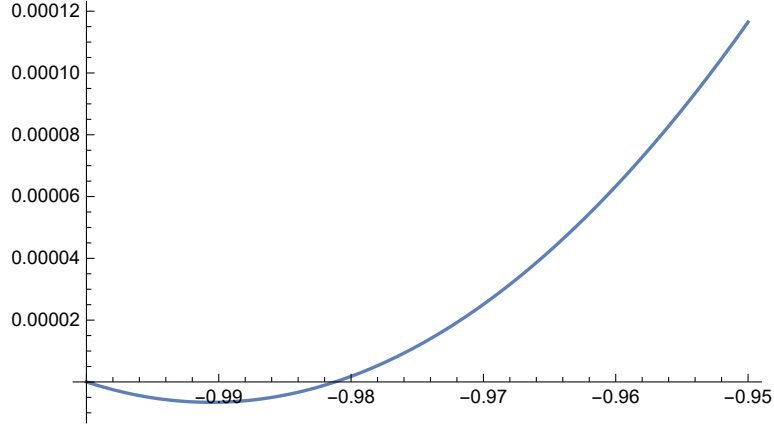}
\caption{The graph of the polynomial $Q_1(x)$ for $b=b_0$}
\label{fig:Q1}
\end{center}
\end{figure}

Letting $N=9,$ we can show the following result by using this strategy
with the aid of computer.

\begin{thm}\label{thm:4/3}
For the logarithmic coefficients $\gamma_n$ of a function $f\in\es,$
the inequality
$$
\sum_{n=1}^\infty \frac{n^2}{n^2+4/3}|\gamma_n|^2
\le B_{2/\sqrt 3}=
\frac{(2\pi/\sqrt{3}) \coth(2\pi/\sqrt 3) -1}{8/3}
=0.98727\cdots
$$
holds, where the inequality is strict unless $f$ is a rotation of the Koebe function.
\end{thm}

\begin{pf}
Let $\varphi(x)=x/(x^2+4/3)$ and $p_n=\varphi(n).$
Since
$$
\varphi''(x)=\frac{54x(x+2)(x-2)}{(3x^2+4)^2}
$$
we find that $\varphi(x)$ is convex on $2\le x<+\infty.$
Note that $p_1-2p_2+p_3=-27/868<0.$
We computed the polynomials $Q_k(x)~(k=1,2,\dots, 9)$ in Theorem \ref{thm:N}
with $N=9$ by using Mathematica as shown in Appendix.
By numerical computations, we found that $Q_k(x)$ has no real roots
for each odd $k$ and that $Q_k(x)$ has only one real root, which is less than $-1,$
for each even $k.$
Thus we confirmed numerically that $Q_k(x)>0$ for $-1<x<1$ and $k=1,\dots, 9.$
We now apply Theorem \ref{thm:N} to get the assertion.
\end{pf}

In a similar way, we can show the following result, which will be used in the next section.
Its proof will also given in Appendix.

\begin{thm}\label{thm:1/20}
Let $\beta=1/20.$
For the logarithmic coefficients $\gamma_n$ of a function $f\in\es,$
the sharp inequality
$$
\sum_{n=1}^\infty \frac{n^3|\gamma_n|^2}{(n+1)^2(n+\beta)}
\le E_{1,\beta}
=\frac{20}{19^2}\left(1-\gamma-\psi\left(\frac{21}{20}\right)\right)-\frac{20}{19}
\left(\frac{\pi^2}6-1\right)
=0.62787\cdots
$$
holds.
\end{thm}

\begin{pf}
Let $\varphi(x)=x^2/[(x+1)^2(x+\beta)]$ with $\beta=1/20$
and $p_n=\varphi(n).$
Then
$$
\varphi''(x)=
\frac{2 \left\{x^2(x-3)(x+1)+(3-6 \beta)  x^2-2 \beta ^2 x+\beta ^2\right\}}%
{(x+1)^4 (\beta +x)^3}
$$
is positive for $x\ge 3.$
In this case, indeed, we have
$$
\lambda_1=-\frac{6985}{630252}
\aand
\lambda_2=\frac{12103}{2025810}
$$
for $\lambda _n=p_n-2p_{n+1}+p_{n+2}.$
We take $N=9$ and compute $Q_k(x)~(k=1,\dots, 9)$ as shown in Appendix.
By numerical computations, as in the previous case,
$Q_k(x)$ has no real roots for each odd $k$ and
$Q_k(x)$ has only one real root, which is less than $-1,$
for each even $k.$
Thus we confirm the assertion
in the same way as the previous theorem.
\end{pf}

\section{Applications}

Our next result is related to a transform $h_f$ of $f\in {\mathcal S}$
introduced by  Danikas and Ruscheweyh \cite{DR99}:
$$
h_f(z):=\int_{0}^{z}\frac{tf'(t)}{f(t)}\,dt.
$$
It was conjectured in \cite{DR99} that the transform $h_f\in {\mathcal S}$
for each $f\in {\mathcal S}$.
This conjecture remains open.
Roth \cite{Roth07} applied his inequality \eqref{eq:Roth} to obtain the sharp
$H^2$ norm estimate of $h_f$ for $f\in\es.$
We now introduce the class
$$
{\mathcal U} = \left \{ f \in {\mathcal A} :\, \left |U_{f}(z)\right | < 1
~\mbox{ for $z\in \D$} \right \},
$$
where
\begin{equation}\label{3-10eq4}
U_{f}(z)=f'(z)\left (\frac{z}{f(z)} \right )^{2}-1, \quad z\in \D.
\end{equation}
It is known that $\U\subset\es.$
See \cite{Aks58} and  also \cite{FP06, OP01, OP07a} and the references therein.
We will say that $f\in \U$ on $|z|<r$ if $f_r(z)=f(rz)/r$ belongs to $\U.$
Several generalizations of the class  $\mathcal{U}$ were investigated in the literature.
Among them, the following result was proved in \cite{OP14}.

\begin{ThmA}[$\text{\cite[Thoerem 4]{OP14}}$]\label{3-13th4}
Let $f\in {\mathcal S}$, $b=|f''(0)|/2!$, and let $H$ be defined by the quotient
\begin{equation}\label{3-10eq7}
H(z)=\frac{z^{2}}{h_f(z)} .
\end{equation}
Then $H\in \mathcal{U}$ on the disk $|z|<r_1(b).$
Here, $r_1(b)\geq r_1(0)\approx 0.557666 $ is the root of the equation
\begin{equation*}\label{eq4a}
\left(\frac{2\pi^2}3-4-\frac{b^2}4\right)r^4(1+r^2)=(1-r^2)^3
\end{equation*}
in $0<r<1$ for $b\in[0,2].$
\end{ThmA}

The proof of this theorem is based on the Roth inequality \eqref{eq:Roth}.
It is almost the optimal choice but there is still room to improve a little as follows.

\begin{thm}\label{thm1}
Let $f\in {\mathcal S}$, $b=|f''(0)|/2!$, and let $H$ be defined by \eqref{3-10eq7}.
Then $H\in \mathcal{U}$ in the disk $|z|<r_2(b)$,
where $r_2(b)\geq r_2(0)\approx 0.558509 $  is the solution of the equation
$$
\log\frac1{1-r^2}+\frac{-r^2+23r^4+18r^6}{(1-r^2)^3}
=\frac{20}{4E_{1,1/20}-5b^2/84}
$$
in $0<r<1$ for $b\in[0,2]$ and $E_{1,1/20}$ is
the constant given in Theorem \ref{thm:1/20}.
\end{thm}

The method of the proof is along the line of \cite{OP14} but based
on Theorem \ref{thm:1/20} instead of the Roth inequality.

\begin{pf}[Proof of Theorem \ref{thm1}]
First we note the expression
$$
h_f(z)=\int_{0}^{z}\left(1+t\left(\log\frac{f(t)}{t}\right)'\right)dt =
z+2\sum_{n=2}^{\infty}\frac{n-1}{n}\gamma _{n-1}z^{n},
$$
where $\gamma _{n}$ $(n\geq 1)$ denote the logarithmic coefficients of
$f\in {\mathcal S}$ defined by \eqref{eq:log}.
We also have
$$
\frac{z}{H(z)}=1+2\sum_{n=1}^{\infty}\frac{n}{n+1}\gamma _{n}z^{n}
$$
and $2|\gamma _1|=|f''(0)|/2=b$.
By the forms of $U_H(z)$ and $H$,  we compute
$$
U_H(z) =-z^{2}\left(\frac{1}{H(z)}-\frac{1}{z} \right )'
= -z\left(\frac{z}{H(z)}\right)'+\frac{z}{H(z)}-1
= -2\sum_{n=2}^{\infty} \frac{(n-1)n}{n+1}\gamma_{n}z^{n}.
$$
Letting $r=|z|<1,$ we estimate with the help of the Cauchy-Schwarz inequality in addition to
Theorem \ref{thm:1/20} as follows:
\begin{eqnarray*}
|U_H(z)|& \leq & 2\sum_{n=2}^{\infty} \frac{(n-1)n}{n+1}|\gamma_{n}|\,|z|^{n}\\
& \leq & 2
\left( \sum_{n=2}^{\infty}\frac{n^3}{(n+1)^2(n+\beta)}|\gamma_{n}|^{2}\right)^{1/2}
\left( \sum_{n=2}^{\infty}\frac{(n-1)^{2}(n+\beta)}{n}r^{2n}\right)^{1/2}\\
& \le & 2\left( E_{1,\beta}-\frac{5|\gamma_{1}|^{2}}{84}\right)^{1/2}
\left(\frac{-r^2+23r^4+18r^6}{20(1-r^2)^3}-\frac1{20}\log(1-r^2)\right)^{1/2}
\end{eqnarray*}
which is less than $1$ whenever,
$$
\log\frac1{1-r^2}+\frac{-r^2+23r^4+18r^6}{(1-r^2)^3}
<\frac{20}{4E_{1,1/20}-5b^2/84}.
$$
Note that the left-hand quantity is increasing from $0$ to $+\infty$
when $r$ moves from $0$ to $1$  so that
the root $r_2(b)$ of the equation in the statement is an increasing function of $b$
on the interval $[0,2].$
\end{pf}

By using Mathematica, we made graphs of the functions $r=r_1(b)$ and $r=r_2(b)$
and a graph of the difference $r_2(b)-r_1(b)$ in Figure 2.

\begin{figure}
\begin{center}
\includegraphics[width=75mm]{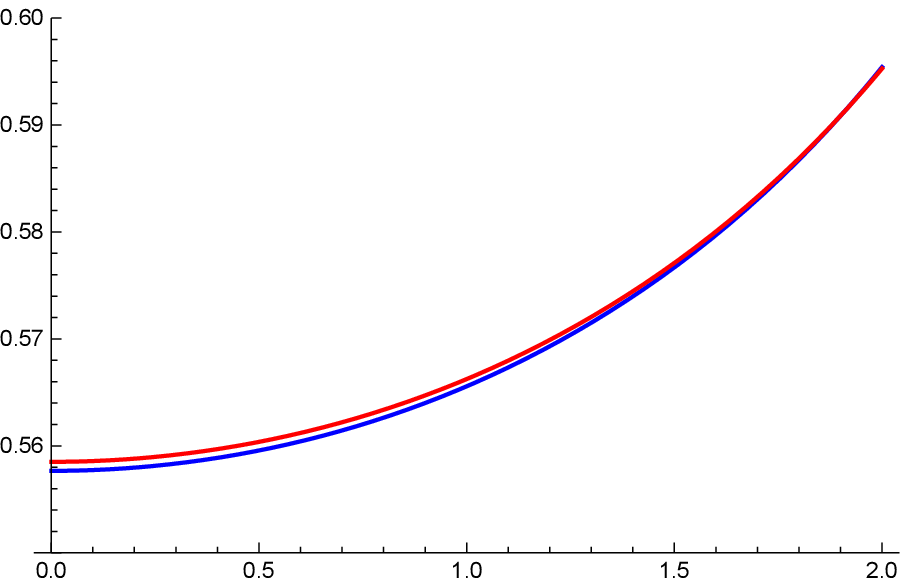}
\includegraphics[width=75mm]{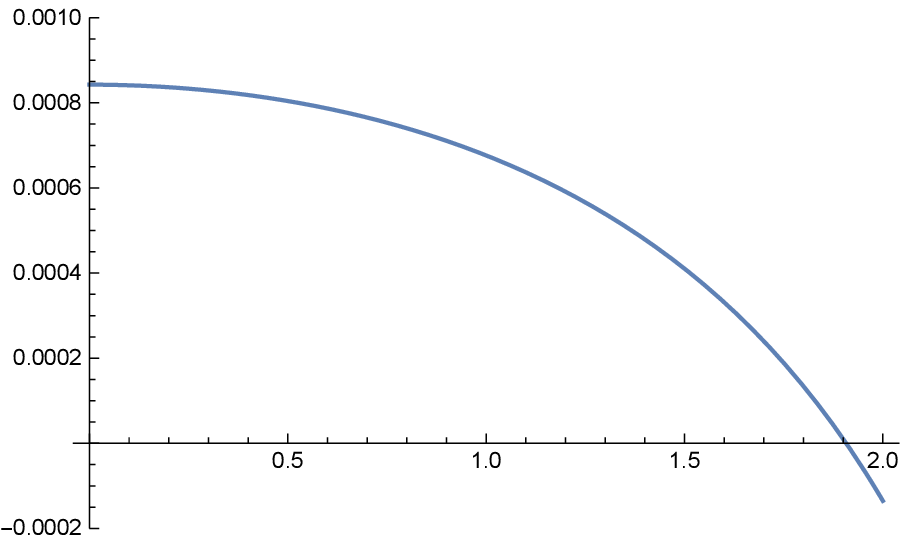}
\caption{Left: the graphs of $r_1(b)$ (blue colored)  and $r_2(b)$ (red colored),
Right: the graph of the difference $r_2(b)-r_1(b)$}
\label{fig:graph}
\end{center}
\end{figure}

In a paper \cite{OPW14}, analytic and geometric properties of the function
$P_f(z)=f(z)/f'(z)$ are studied for $f\in\es.$
Let us look at the following result in the paper.

\begin{ThmB}[$\text{\cite[Theorem 3.16]{OPW14}}$]
Let $f\in\es$ and $b=|f''(0)|/2!.$
Then $P_f\in\U$ on the disk $|z|<r_3(b),$ where $r=r_3(b)\ge r_3(0)\approx
0.360794$ is the root of the equation
$$
\left(\frac{2\pi^2}3-4-\frac{b^2}4\right)r^4(r^6-5r^4+19r^2+9)=(1-r^2)^5
$$
in $0<r<1$ for $b\in[0,2].$
\end{ThmB}

Their proof relied also on the Roth inequality \eqref{eq:Roth}.
Here, we replace it by Theorem \ref{thm:4/3}.

\begin{thm}\label{thm2}
Let $f\in\es$ and $b=|f''(0)|/2!.$
Then $P_f\in\U$ on the disk $|z|<r_4(b).$
Here, $r=r_4(b)$ is the solution of the equation
$$
\left(4B_{2/\sqrt 3}-\frac{3b^2}7\right)r^4(r^6+2r^4+11r^2+4)=\frac34(1-r^2)^5
$$
in $0<r<1$ for $b\in[0,2]$ and $B_{2/\sqrt 3}$ is the constant given in Theorem \ref{thm:4/3}.
The function $r_4(b)$ is increasing in $0\le b\le 2$ and $r_4(b)\ge r_4(0)\approx 0.362012.$
\end{thm}

\begin{pf}
Let $F=P_f.$
Since
$$
\frac{zf'(z)}{f(z)}=1+2\sum_{n=1}^\infty n\gamma_nz^n,
$$
we obtain the expressions
\begin{align*}
U_F(z)&=F'(z)\left (\frac{z}{F(z)} \right )^{2}-1
=\frac{zf'(z)}{f(z)}-z\left (\frac{zf'(z)}{F(z)} \right )'-1 \\
&=-2\sum_{n=1}^\infty n(n-1)\gamma_n z^n.
\end{align*}
Hence, as in the proof of Theorem \ref{thm1}, we estimate
\begin{eqnarray*}
|U_H(z)|& \leq & 2\sum_{n=2}^{\infty} (n-1)n|\gamma_{n}|\,|z|^{n}\\
& \leq & 2
\left( \sum_{n=2}^{\infty}\frac{n^2}{n^2+4/3}|\gamma_{n}|^{2}\right)^{1/2}
\left( \sum_{n=2}^{\infty}(n-1)^{2}(n^2+4/3)r^{2n}\right)^{1/2}\\
& \le & 2\left( B_{2/\sqrt 3}-\frac{3|\gamma_{1}|^{2}}{7}\right)^{1/2}
\left(\frac{4r^4(r^6+2r^4+11r^2+4)}{3(1-r^2)^5}\right)^{1/2}.
\end{eqnarray*}
We now see that $|U_F(z)|<1$ as long as
$$
\frac{4r^4(r^6+2r^4+11r^2+4)}{3(1-r^2)^5}
<\frac1{4B_{2/\sqrt 3}-3b^2/7}.
$$
Now the assertion follows as before.
\end{pf}

In Figure 3, we exhibit the graphs of $r_3(b)$ and $r_4(b).$

\begin{figure}
\begin{center}
\includegraphics[width=80mm]{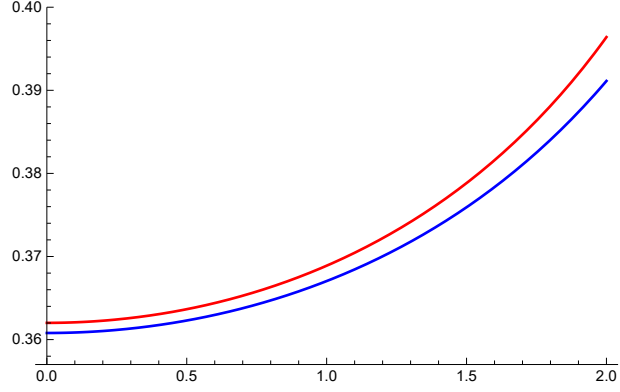}
\caption{The graphs of $r_3(b)$ (blue colored) and $r_4(b)$ (red colored)}
\label{fig:graph3}
\end{center}
\end{figure}

\newpage
\section{Appendix}

The polynomials $Q_k(x)~(k=1,2,\dots, 9)$ used in the proof of
Theorem \ref{thm:4/3} are presented below.
We note that by using a suitable command of Mathematica or similar software,
we can find all the roots of the following polynomials numerically.
In this way, we can check that $Q_k$ has no roots on the interval $(-1,1)$ so that
$Q_k(x)>0$ for $(-1,1).$
\begin{eqnarray*}
&Q_1(x)&=\frac{1136025 x^8}{3570176}+\frac{387585 x^7}{575456}+\frac{285789141 x^6}{943419008}-\frac{103110975 x^5}{825491632}+\frac{34505962335 x^4}{1043421422848}\\
&&+\frac{380568045735 x^3}{1695559812128}+\frac{46758786465915 x^2}{210249416703872}+\frac{2558807811009 x}{13140588543992}+\frac{77049161884395}{840997666815488}, \\
&Q_2(x)&=\frac{103275 x^7}{223136}+\frac{3957741 x^6}{3123904}+\frac{9202167 x^5}{7370461}+\frac{108476415 x^4}{173787712}+\frac{7350044085 x^3}{18632525408}\\
&&+\frac{1429990103205 x^2}{3391119624256}+\frac{9108847966527 x}{26281177087984}+\frac{19344079210563}{105124708351936}, \\
&Q_3(x)&=\frac{240975 x^6}{446272}+\frac{2618811 x^5}{1561952}+\frac{984379149 x^4}{471709504}+\frac{1223639757 x^3}{825491632}+\frac{222176774097 x^2}{260855355712}\\
&&+\frac{821182609953 x}{1695559812128}+\frac{20630684258217}{105124708351936},\\
&Q_4(x)&=\frac{722925 x^5}{1450384}+\frac{32587515 x^4}{20305376}+\frac{805347441 x^3}{383263972}+\frac{16309308609 x^2}{10731391216}+\frac{607240980387 x}{847779906064}\\
&&+\frac{357105897585}{1695559812128}, \\
&Q_5(x)&=\frac{516375 x^4}{1450384}+\frac{696195 x^3}{634543}+\frac{1015598331 x^2}{766527944}+\frac{1064477133 x}{1341423902}+\frac{188723150883}{847779906064},\\
&Q_6(x)&=\frac{34425 x^3}{181298}+\frac{9077265 x^2}{17767204}+\frac{326148282 x}{670711951}+\frac{463541805}{2682847804},\\
&Q_7(x)&=\frac{103275 x^2}{1450384}+\frac{5178573 x}{35534408}+\frac{857191005}{10731391216},\\
&Q_8(x)&=\frac{6075 x}{362596}+\frac{684531}{35534408}, \\
&Q_9(x)&=\frac{675}{362596}.
\end{eqnarray*}

\newpage

The polynomials $Q_k(x)~(k=1,2,\dots, 9)$ used in the proof of
Theorem \ref{thm:1/20} are presented below.

\begin{eqnarray*}
&Q_1(x)&=\frac{62987827 x^8}{341496320}+\frac{330128990251 x^7}{860326246080}+\frac{158362220519819 x^6}{1016659815367680}\\
&&-\frac{10529204214766063 x^5}{107935383731535360}-\frac{530240051345429 x^4}{93380566898691072}+\frac{541351055064272599 x^3}{5450736878442535680}\\
&&+\frac{1284708860080692110137 x^2}{14629777781739765765120}+\frac{1777709878968276897929 x}{27264585865969563471360}\\
&&+\frac{1284524775560504080639}{43889333345219297295360}, \\
&Q_2(x)&=\frac{62987827 x^7}{234778720}+\frac{5390982010279 x^6}{7435676841120}+\frac{1932000965815939 x^5}{2795814492261120}\\
&&+\frac{110402609517184133 x^4}{356186766314066688}+\frac{438453616394854873 x^3}{2569633099837195392}+\frac{16605230695470458161 x^2}{89937158494301838720}\\
&&+\frac{1581160573987373624339 x}{10972333336304824323840}+\frac{30873002384858864057389}{449865666788497797277440}, \\
&Q_3(x)&=\frac{440914789 x^6}{1408672320}+\frac{2383515042599 x^5}{2478558947040}+\frac{3620224806455 x^4}{3089297781504}\\
&&+\frac{70911808727134489 x^3}{89046691578516672}+\frac{852605666350380625 x^2}{1998603522095596416}+\frac{10225388053815434921 x}{44968579247150919360}\\
&&+\frac{942885972184285561567}{10972333336304824323840}, \\
&Q_4(x)&=\frac{440914789 x^5}{1526061680}+\frac{8910108245315 x^4}{9666379893456}+\frac{98087517932455 x^3}{82603609998624}\\
&&+\frac{20630163482121289 x^2}{24735192105143520}+\frac{104078041421486463 x}{277583822513277280}+\frac{2331963676119479803}{22484289623575459680}, \\
&Q_5(x)&=\frac{62987827 x^4}{305212336}+\frac{1524238766119 x^3}{2416594973364}+\frac{105039968099029 x^2}{139790724613056}\\
&&+\frac{5760631998672587 x}{13095101702723040}+\frac{486393597309120707}{4088052658831901760}, \\
&Q_6(x)&=\frac{62987827 x^3}{572273130}+\frac{318698320009 x^2}{1084369539330}+\frac{4395382293406319 x}{15901194924735120}+\frac{10724868335939851}{111308364473145840}, \\
&Q_7(x)&=\frac{62987827 x^2}{1526061680}+\frac{14188797597869 x}{169161648135480}+\frac{120529899535861}{2650199154122520}, \\
&Q_8(x)&=\frac{62987827 x}{6485762140}+\frac{938138865611}{84580824067740}, \\
&Q_9(x)&=\frac{62987827}{58371859260}.
\end{eqnarray*}

\def\cprime{$'$} \def\cprime{$'$} \def\cprime{$'$}
\providecommand{\bysame}{\leavevmode\hbox to3em{\hrulefill}\thinspace}
\providecommand{\MR}{\relax\ifhmode\unskip\space\fi MR }
\providecommand{\MRhref}[2]{%
  \href{http://www.ams.org/mathscinet-getitem?mr=#1}{#2}
}
\providecommand{\href}[2]{#2}

\end{document}